\documentclass[12pt]{amsart}
\usepackage{pstricks}
\usepackage{amsfonts}
\usepackage{amssymb}
\usepackage{amsthm}
\usepackage{latexsym,amsmath}
\usepackage{graphicx}
\usepackage{stmaryrd}

\newtheorem{theorem}{Theorem}[section]
\theoremstyle{definition}
\newtheorem{Definition}[theorem]{Definition}

\newenvironment{theorem*}[1]{\medskip
                            \noindent
                            {\bf Theorem #1. }\ %
                            \begingroup \sl}
                            {\endgroup\medskip}

\title
 {Which axioms of set theory are intrinsically justified?}

\author[$\mathrm{M^{\lowercase{c}}Callum}$ ]{\textbf{Rupert} $\mathbf{M^{\lowercase{c}}Callum}$ }

\begin{document}

\begin{abstract}

We recently formulated a new large-cardinal axiom of strength intermediate between a totally indescribable cardinal and an $\omega$-Erd\H{o}s cardinal, positing the existence of what we called an ``extremely reflective cardinal", and we showed that the property of being extremely reflective was in fact equivalent to the property of being remarkable, and we sought to argue that this axiom should be seen as intrinsically justified. This built on related earlier work in which the notion of an $\alpha$-reflective cardinal was formulated. Then Welch and Roberts put forward a family of reflection principles, Welch's principle implying the existence of a proper class of Shelah cardinals and provably consistent relative to a superstrong cardinal, and Roberts' principle implying the existence of a proper class of 1-extendible cardinals and provably consistent relative to a 2-extendible cardinal. Roberts tentatively argued that his principle should be seen as intrinsically justified (at least on the assumption that a weaker form of reflection involving reflection of second-order formulas with a second-order parameter should be seen as intrinsically justified). This work overlapped with previous work of Victoria Marshall's on reflection principles. We analyze the relationship between reflection principles equivalent to those studied in my earlier work and stronger but similar reflection principles which are natural extensions of those of Welch and Roberts. We also show how a natural strengthening of Roberts' reflection principle yields the existence of supercompact cardinals, and in the process solve a question which Marshall left open, of whether her theory $B_0(V_0)$ is strong enough to imply the existence of supercompact cardinals. We also manage to resolve negatively her question of whether her theory $B_1(V_0)$ implies the existence of a huge cardinal.

\end{abstract}

\maketitle

\section{Introduction}

In speaking of justifications for candidates for new axioms for set theory to be added to the ZFC axioms, philosophers of set theory draw a distinction between ``intrinsic" and ``extrinsic" justifications. This distinction has its origins in remarks made by Kurt G\"odel in his famous essay about the continuum problem \cite{Godel1964}. There he wrote `First of all, the axioms of set theory by no means form a system closed in itself, but, quite on the contrary, the very concept of set on which they are based suggests their extension by new axioms which assert the existence of still further iterations of the operation ``set of". These axioms can be formulated also as propositions asserting the existence of very great cardinal numbers (i.e., of sets having these cardinal numbers). The simplest of these strong ``axioms of infinity" asserts the existence of inaccessible numbers (in the weaker or stronger sense) $> \aleph_{0}$. The latter axiom, roughly speaking, means nothing else but that the totality of sets obtainable by use of the procedures of formation of sets expressed in the other axioms forms again a set (and, therefore, a new basis for further applications of these procedures). Other axioms of infinity have first been formulated by P. Mahlo. These axioms show clearly, not only that the axiomatic system of set theory as used today is incomplete, but also that it can be supplemented without arbitrariness by new axioms which only unfold the content of the concept of set explained above.' In a footnote he wrote `... some propositions have been formulated which, if consistent, are extremely strong axioms of infinity of an entirely new kind... That these axioms are implied by the general concept of set in the same sense as Mahlo's has not been made clear yet.' Here he was referring to axioms asserting the existence of weakly compact, measurable, or strongly compact cardinals.

\bigskip

Thus, he held the view that axioms asserting the existence of inaccessible or Mahlo cardinals were ``intrinsically justified" in the sense of merely unfolding the content of the conception of set, whereas it was not clear whether this was the case for axioms asserting the existence of weakly compact, measurable, or strongly compact cardinals, and at the time of his writing the essay the only known justifications for these axioms were ``extrinsic". An extrinsic justification for a candidate for a new axiom is based on something other than the idea that the axiom merely unfolds the content of the conception of set, such as for example the axiom having desirable or verifiable consequences. G\"odel made the following remarks about this kind of justification: `Secondly, however, even disregarding the intrinsic necessity of some new axiom, and even in case it has no intrinsic necessity at all, a probable decision about its truth is possible also in another way, namely, inductively by studying its ``success". Success here means fruitfulness in consequences, in particular in ``verifiable" consequences, i.e., consequences demonstrable without the new axiom, whose proofs with the help of the new axiom, however, are considerably simpler and easier to discover, and make it possible to contract into one proof many different proofs... There might exist axioms so abundant in their verifiable consequences, shedding so much light upon a whole field, and yielding such powerful methods for solving problems (and even solving them constructively, as far as that is possible) that, no matter whether or not they are intrinsically necessary, they would have to be accepted at least in the same sense as any well-established physical theory.'

\bigskip

These remarks were the origin of the distinction between ``intrinsic" and ``extrinsic" justifications which is still discussed by philosophers of set theory today.

\bigskip

In a recent paper \cite{McCallum2017} I introduced the notion of an ``extremely reflective cardinal", and proved that the property of being extremely reflective was in fact equivalent to the property of being remarkable. For that reason I should say ``remarkable cardinal" from now on rather than ``extremely reflective cardinal". I made a few remarks about intrinsic justifications for large-cardinal axioms in general, influenced by Tait's ideas in \cite{Tait2005a}, supporting G\"odel's view that axioms asserting the existence of inaccessible and Mahlo cardinals are indeed intrinsically justified, in the sense described above. Then I used these ideas to motivate the view that belief in the existence of a remarkable cardinal was also intrinsically justified. This work built on related work on the topic of intrinsically justified reflection principles by Tait \cite{Tait2005a} and Koellner \cite{Koellner2009}. The challenge of formulating an intrinsically justified reflection principle with at least as much consistency strength as an $\omega$-Erd\H{o}s cardinal was made by Peter Koellner in \cite{Koellner2009}. Roberts recently wrote a very interesting paper \cite{Roberts2017} in response to this challenge, building on previous similar work of Welch in \cite{Welch2014} but avoiding explicit reference to elementary embeddings.

\bigskip

To explain the reflection principle which Roberts formulates in his paper, let us begin by explaining the reflection principle that he calls $\textsf{R}_{2}$. This is an axiom schema in the second-order language of set theory. For each formula $\phi(x_{1}, x_{2}, \ldots x_{m}, X_{1}, X_{2}, \ldots X_{n})$ in the second-order language of set theory, there is an axiom asserting that if $\phi$ holds, then there exists an ordinal $\alpha$ such that $x_{1}, x_{2}, \ldots x_{m} \in V_{\alpha}$, and a ``set-sized" family of classes which contains the classes $X_{1}, X_{2}, \ldots X_{n}$, which is itself coded for by a single class, and which is standard for $V_{\alpha}$ in the sense that every subset $X\subseteq V_{\alpha}$ is such that some class in the family has intersection with $V_{\alpha}$ equal to $X$, such that the formula $\phi$ still holds when the first-order variables are relativised to $V_{\alpha}$ and the second-order variables are relativised to the set-sized family of classes. This completes the description of the axiom schema $\textsf{R}_{2}$. Then Roberts extends the axiom schema as follows. He extends the underlying language so as to include a satisfaction predicate for the second-order language of set theory, and then he extends the axiom schema so as to also include an axiom of the kind described for every formula in this extended language, calling this new axiom schema $\textsf{R}_{S}$. Then he denotes by $\textsf{ZFC2}_{S}$ the result of extending $\textsf{ZFC2}$ -- being the same as ZFC except for having Separation and Replacement as single second-order axioms and also having an axiom schema of class comprehension for every formula in the second-order language of set theory -- by adding the usual Tarskian axioms for the satisfaction predicate and extending the class comprehension axiom schema to include axioms involving formulas in the extended language. Then he proceeds to investigate the theory $\textsf{ZFC2}_{S}+\textsf{R}_{S}$. This completes the description of the reflection principle which Roberts considers. He shows that the theory $\textsf{ZFC2}_{S}+\textsf{R}_{S}$ proves the existence of a proper class of 1-extendible cardinals and is consistent relative to a 2-extendible cardinal. The reflection principle suggested by Welch in \cite{Welch2014} is similar but slightly weaker, and makes use of an elementary embedding rather than a ``reflecting structure".

\bigskip

In Section 2 I will describe the ideas about intrinsic justification which originally led me to think that the existence of remarkable cardinals was intrinsically justified. Then in Section 3 I will analyze the relationship between the reflection principles of \cite{McCallum2013} and \cite{McCallum2017} with those of \cite{Welch2014}, \cite{Roberts2017} and \cite{Marshall89}, and consider the philosophical question whether there can be any principled grounds for accepting the former but not the latter. I will also discuss how a natural extension of the reflection principle put forward by Roberts is equivalent to the existence of a supercompact cardinal. In Section 4 I make some concluding remarks.

\section{Intrinsic justifications for small large cardinals}

It is instructive to begin by considering a challenge posed by Schindler in \cite{Schindler1994} to the idea that any large cardinals are intrinsically justified. In that paper Schindler considers a class theory $BL_{1}$ which he claims to embody the weakest reflection principle schema that leads to the existence of large cardinals. He shows that it proves $\Delta^{1,NBG}_{1}$-comprehension but that this comprehension schema is false for the predicative classes. He concludes that we then face the dilemma between accepting the existence of non-predicative classes (for which he believes there is no good justification), or concluding that there is not sufficient justification for believing in the existence of large cardinals, at least in $V$ (which he finds to be an unsatisfying view). Let us now see how we can respond to this argument based on ideas of Tait in \cite{Tait2005a}.

\bigskip

In \cite{Tait2005a} Tait begins by quoting Cantor in \cite{Cantor1883} as saying, with reference to the totality $\Omega$ of transfinite ordinal numbers, `If the initial segment $\Sigma$ of $\Omega$ is a set, then it has a least strict upper bound $S(\Sigma)\in\Omega$.' It is obviously not consistent for every initial segment of $\Omega$ to be a set, since $\Omega$ itself cannot be a set. So some criterion must be given for when an initial segment is a set. Later Cantor attempted to solve this problem by invoking the notion of a `consistent multiplicity', but as Tait says `that is just naming the problem, not solving it.'  Tait proposes a modified version of Cantor's principle which Peter Koellner later in \cite{Koellner2009} called the Relativised Cantorian Principle, which says `If the initial segment $\Sigma$ of $\Omega_{C}$ satisfies the condition $C$, then it has a least strict upper bound $S(\Sigma)\in\Omega_{C}$.' With an appropriate choice of the condition $C$, which Tait calls an `existence condition', we can in this way determine an initial segment $\Omega_{C}$ of the totality of transfinite ordinal numbers $\Omega$. Tait writes `In this way, we are led to a hierarchy of more and more inclusive existence conditions, each of which can replace the condition ``\textit{is a set}" in Cantor's definition and so yields an initial segment of the transfinite numbers, but none of which yields ``all the numbers".' He calls this the `bottom-up' conception.

\bigskip

We can articulate the idea that this is part of the iterative conception of set as follows. We view ourselves as building up levels $V_{\alpha}$ of the universe from below. If, at any given stage of the process, we have built up a level $V_{\alpha}$ which satisfies a certain condition which we take to be a ``warrant" for concluding that it cannot be all of $V$, then we take ourselves to be justified in going one level further. We also take ourselves to be justified in taking limits at limit ordinals. Furthermore, given any circumscribed set of processes for building up new levels of the universe in this way, we take ourselves to be justified in positing a closure point for all these processes which is still only a set. In this way we can motivate the proposition that a level satisfying a certain amount of reflection should exist, (that is, properties of the level are reflected to a lower level), by positing that when a level fails to satisfy that form of reflection then that is warrant for going one level further.

\bigskip

This can be used to motivate the existence of a level satisfying a first-order reflection schema, which is equivalent, given the axioms of ZF other than Infinity and Replacement, to all of ZF. But we can also take ourselves to be justified in positing reflection schemas for formulas in higher-order languages, because at each stage of the process of building up the levels from below, we know that the level we have built up so far is only a set and so we are justified in speaking of properties of this level which can only be formulated in higher-order languages. Of course, this knowledge that the level we have built up so far is only a set cannot by itself be taken to be warrant for going one level further, otherwise the principles we have enunciated so far would lead to inconsistency. Rather, there must be some condition such that, if the level built up so far satisfies it, then that is a warrant for concluding that it cannot be all of $V$ and so going one level further. This can be based on the idea that $V$ should not be ``describable" and so any witness to the ``describability" of this level can be taken as a warrant for going one level further. This line of thought can be used to motivate higher-order reflection principles (or at least levels of the universe satisfying higher-order reflection principles) without any need to be committed to the existence of non-predicative proper classes. That is, these ideas of Tait provide us with a way of responding to Schindler's critique of the view that belief in large cardinals in justified.

\bigskip

If we accept that higher-order reflection principles can be motivated in this way, the question then arises which reflection principles should be accepted. One naive idea might be as follows. Let $\mathcal{L}_{\in}^{\omega}$ be the $\omega$th-order language of set theory, and suppose we are considering a formula $\phi$ with free variables $x_1, x_2, \ldots x_{m_{1}}, X^{(2)}_{1}, X^{(2)}_{2}, \ldots X^{(2)}_{m_{2}}, X^{(3)}_{1}, X^{(3)}_{2}, \ldots \newline X^{(3)}_{m_{3}},  \ldots X^{(k)}_{1}, X^{(k)}_{2}, \ldots X^{(k)}_{m_{k}},$, where the variables denoted by lower-case letters are first-order variables, and for variables denoted by upper-case letters superscripts indicate the order of the variables, and an $m$th-order variable ranges over $\mathcal{P}^{m-1}(D)$ where $D$ is the domain of discourse. Suppose that the domain of discourse $D$ is equal to $V_{\kappa}$ for some ordinal $\kappa$. If $X^{(2)}\in\mathcal{P}(V_{\kappa})$, let $X^{(2),\alpha}:=X^{(2)}\cap V_{\alpha}$ for all ordinals $\alpha<\kappa$, and then inductively define $X^{(k),\alpha}=\{X^{(k-1),\alpha}\mid X^{(k-1)}\in X^{(k)}\}$ for all ordinals $\alpha<\kappa$ and all $X^{(k)}\in\mathcal{P}^{k-1}(V_{\kappa})$. Then we can posit that there should be a level $V_{\kappa}$ with the property that if the formula $\phi$ is true in $V_{\kappa}$ for a certain assignment of values to the variables, then there should be some ordinal $\alpha<\kappa$ such that the relativisation of $\phi$ to $V_{\alpha}$ holds for an assignment of values to the variables such that if a variable has $X^{(j)}$ assigned to it under the first assignment then it will have $X^{(j),\alpha}$ as defined above assigned to it under the second assignment. By the relativisation of $\phi$ to $V_{\alpha}$ we mean the result of relativising the first-order quantifiers to $V_{\alpha}$, the second-order quantifiers to $V_{\alpha+1}$, and so on. This is a naive form of reflection which one might naturally posit.

\bigskip

If one makes the restriction that the free variables should be no higher than second-order then this is consistent (relative to the consistency of a totally indescribable cardinal), but if one allows parameters of third order or higher then this is inconsistent, as was first observed by Reinhardt in \cite{Reinhardt74}. This can be seen as follows, using an argument presented by Tait in \cite{Tait2005a}. Let $\phi$ be a formula with exactly one free third-order variable $X^{(3)}$ asserting that every element of $X^{(3)}$ is a bounded set of ordinals. Then if we assign to the variable $X^{(3)}$ the set $\{\{\alpha\mid\alpha<\beta\}\mid\beta<\kappa\}$, then $\phi(X^{(3)})$ will be true in $V_{\kappa}$ but we will not have $\phi^{\alpha}(X^{(3),\alpha})$ for any $\alpha<\kappa$, where $\phi^{\alpha}$ denotes the relativisation of $\phi$ to $V_{\alpha}$. Thus the reflection principle described is inconsistent if we allow parameters of third order or higher. So if one wants to accept the reflection principle in the case of parameters no higher than second order, one must identify principled reasons why one should not also accept it for parameters of third order or higher.

\bigskip

Let us consider what Tait writes in \cite{Tait2005a} about this question.

\bigskip

`One plausible way to think about the difference between reflecting $\phi(A)$ when $A$ is second-order and when it is of higher-order is that, in the former case, reflection is asserting that, if $\phi(A)$ holds in the structure $\langle R(\kappa),\in,A \rangle$, then it holds in the substructure $\langle R(\beta),\in, A^{\beta} \rangle$ for some $\beta<\kappa \ldots$ But, when $A$ is higher-order, say of third-order this is no longer so. Now we are considering the structure $\langle R(\kappa), R(\kappa+1), \newline \in, A \rangle$ and $\langle R(\beta), R(\beta+1), \in, A^{\beta} \rangle$. But, the latter is not a substructure of the former, that is the `inclusion map' of the latter structure into the former is no longer single-valued: for subclasses $X$ and $Y$ of $R(\kappa), X\neq Y$ does not imply $X^{\beta}\neq Y^{\beta}$. Likewise for $X\in R(\kappa+1), X\notin A$ does not imply $X^{\beta}\notin A^{\beta}$. For this reason, the formulas that we can expect to be preserved in passing from the former structure to the latter must be suitably restricted and, in particular, should not contain the relation $\notin$ between second- and third-order objects or the relation $\neq$ between second-order objects.'

\bigskip

He then uses these ideas to motivate the following family of reflection principles.

\begin{Definition} A formula in the $\omega$th-order language of set theory is positive iff it is built up by means of the operations $\vee$, $\wedge$, $\forall$, $\exists$ from atoms of the form $x=y$, $x \neq y$, $x \in y$, $x \notin y$, $x \in Y^{(2)}$, $x \notin Y^{(2)}$ and $X^{(m)} = X'^{(m)}$ and
$X^{(m)} \in Y^{(m+1)}$, where $m \geq 2$. \end{Definition}

\begin{Definition} Suppose that $A$ is a variable of order at most 2. Define $A_{n\times}=\{\langle n,x\rangle\mid x\in A\}$ and $A_{/n}=\{x\mid\langle n,x\rangle\in A\}$. In the case where $B$ is a variable of order greater than 2, define, using induction on the order of $B$, $B_{n\times}:=\{A_{n\times}\mid A\in B\}$ and $B_{/n}:=\{A_{/n}\mid A\in B\}$. Given an $n$-tuple of variables of order $m$, $(X_{0}^{(m)}, X_{1}^{(m)}, \ldots X_{n-1}^{(m)})$, call the operator $S$ which sends this $n$-tuple to $S(X_{0}^{(m)}, X_{1}^{(m)}, \ldots X_{n-1}^{(m)}):=\bigcup_{i=0}^{n-1} (X_{i}^{(m)})_{i\times}$ a contracting operator, and call an operator of the form $X \mapsto X_{/n}$ a de-contracting operator. Call a formula $\phi$ in the $\omega$th-order language of set theory positive in the extended sense if it is equivalent to a formula, which has no quantified variables of order higher than first-order, and is positive, except that it may include contracting and de-contracting operators. \end{Definition}

\begin{Definition} For $0<n<\omega$, $\Gamma^{(2)}_{n}$ is the class of formulas

$$ \forall X_{1}^{(2)} \exists Y_{1}^{(k_{1})} \cdots \forall X_{n}^{(2)} \exists Y_{n}^{(k_{n})} \phi(X_{1}^{(2)}, Y_{1}^{(k_{1})}, \ldots, X_{n}^{(2)}, Y_{n}^{(k_{n})}, A^{(l_{1})}, \ldots A^{(l_{n'})}) $$

where $\phi$ is positive in the extended sense and $k_{1}, \ldots k_{n}, l_{1}, \ldots l_{n'}$ are natural numbers which are greater than or equal to 2. \end{Definition}

\begin{Definition} We say that $V_{\kappa}$ satisfies $\Gamma^{(2)}_{n}$-reflection if, for all $\phi\in\Gamma^{(2)}_{n}$, if \newline $V_{\kappa}\models\phi(A^{(m_{1})},A^{(m_{2})},\ldots A^{(m_{p})})$ then $V_{\kappa}\models\phi^{\delta}(A^{(m_{1}),\delta},A^{(m_{2}),\delta},\ldots A^{(m_{p}),\delta})$ for some $\delta<\kappa$, where $\phi^{\delta}$ denotes the relativisation of $\phi$ in the manner described previously. \end{Definition}

Peter Koellner established in \cite{Koellner2009} that these reflection principles are consistent relative to an $\omega$-Erd\H{o}s cardinal. In \cite{Tait2005a} Tait proposes to define $\Gamma^{(m)}_{n}$ in the same way as the class of formulas $\Gamma^{(2)}_{n}$, except that universal quantifiers of order $\leq m$ are permitted. Koellner shows in \cite{Koellner2009} that this form of reflection is inconsistent when $m>2$.

\bigskip

Now I shall describe a family of reflection principles which are consistent relative to an $\omega$-Erd\H{o}s cardinal and which yield all of the reflection principles of Tait considered so far not known to be inconsistent. I shall describe a naive version of it first which is inconsistent and then show how to modify it so as to make it consistent relative to an $\omega$-Erd\H{o}s cardinal.

\bigskip

Consider a formula $\phi$ in the $\omega$th-order language of set theory with free variables $x_{1}, x_{2}, \ldots x_{m_{1}}, X^{(2)}_{1}, X^{(2)}_{2}, \ldots X^{(2)}_{m_{2}}, \ldots X^{(k)}_{1}, X^{(k)}_{2}, \ldots X^{(k)}_{m_{k}}$ and suppose that $\phi$ holds in $V_{\kappa}$ for a certain assignment of values to the free variables. The first idea we try is simply to avoid re-interpreting the higher-order variables when we reflect downwards. That is, we posit that there should be an $\alpha<\kappa$ with $x_{1}, x_{2}, \ldots x_{m_{1}} \in V_{\alpha}$, such that $\phi$ holds with all of the free variables interpreted the same way as before, first-order variables ranging over $V_{\alpha}$ and $k$th-order variables ranging over $\mathcal{P}^{k-1}(V_{\kappa})$ for $k\geq 2$. This is a natural form of reflection to consider as an alterative to the first family of reflection principles we considered that was shown to be inconsistent. However, this form of reflection is also unfortunately inconsistent. To see this, consider the sentence $\phi$ given by $\forall X \forall Y (X \neq Y \implies \exists x ((x \in X \wedge x \notin Y) \vee (x \notin X \wedge x \in Y)))$, where $X$ and $Y$ are second-order variables. This formula clearly gives a counter-example to our reflection principle for any level $V_{\kappa}$.

\bigskip

We can diagnose the root of this problem as arising from the fact that we introduced an existential quantifier for a first-order variable within the scope of a higher-order quantifier. This means that when we ``Skolemize" the formula the existentially quantified first-order variable becomes a Skolem function of a tuple of variables including higher-order variables, and now the fact that the number of possible values for a higher-order variable is large may mean that we cannot have closure under the Skolem functions for any reflecting structure $(V_{\alpha}, V_{\kappa+1}, V_{\kappa+2}, \ldots)$ for any $\alpha<\kappa$. If we modify the reflection principle so as only to apply to formulas in which unbounded existential quantifiers for first-order variables do not appear within the scope of a higher-order quantifier (but existential quantifiers bound by a first-order variable are allowed), then the resulting reflection principle becomes consistent relative to an $\omega$-Erd\H{o}s cardinal. In fact a level $V_{\kappa}$ satisfies this form of reflection if and only if $\kappa$ is $\omega$-refelctive in the sense defined in \cite{McCallum2013}, and when this is so then $V_{\kappa}$ satisfies $\Gamma^{(2)}_{n}$-reflection for all $n \in\omega$.

\bigskip

It is natural to posit this as the correct modification of the inconsistent reflection principle considered earlier. An $\omega$-reflective cardinal is totally ineffable and every totally ineffable cardinal is a stationary limit of totally indescribable cardinals, so we obtain a justification for the form of reflection considered earlier in the case of parameters no higher than second order. As discussed in \cite{McCallum2017} a natural development of this line of thought motivates the idea that a remarkable cardinal is also intrinsically justified.

\bigskip

This completes our discussion of intrinsic justifications for small large cardinals. All cardinals considered so far are indeed small in the sense of being consistent relative to an $\omega$-Erd\H{o}s cardinal. The time has now come to consider the reflection principles of \cite{Roberts2017} and \cite{Marshall89}.

\section{The Reflection Principles of Welch, Roberts and Marshall}

Clearly the reflection principle that I described in the previous section is quite similar to Roberts' reflection principle. If we accept the former as intrinsically justified, then are there any principled reasons to entertain doubts about the latter being intrinsically justified? One possible consideration is the following. If $V_{\kappa}$ is a level of the universe which fails to satisfy the reflection principle that I described in the previous section, as witnessed by some formula $\phi$ with parameters $x_{1}, x_{2}, \ldots x_{m_{1}}, X^{(2)}_{1}, X^{(2)}_{2}, \ldots X^{(2)}_{m_{2}}, \ldots$ $X^{(k)}_{1}, X^{(k)}_{2}, \ldots  X^{(k)}_{m_{k}}$, then the structure $(V_{\kappa}, V_{\kappa+1}, V_{\kappa+2}, \ldots )$ is describable relative to these parameters as the first element in a transfinite well-ordered sequence of structures satisfying the formula $\phi$ in question relative to the parameters in question, where the well-ordering is itself independent of any parameters. In this situation there is canonical well-ordering of the family of structures in question. However, when considering Roberts' principle, there is no canonical well-ordering of the family of possible reflecting structures, and so one may resist the idea that failure of the reflection principle at a particular level shows that that level is ``describable" relative to the parameters. Perhaps this might constitute a principled reason for accepting the reflection principle described in the previous section as intrinsically justified, but refusing to accept Roberts' reflection principle as intrinsically justified.

\bigskip

A possible response is the following. We could take some well-ordering of the family of possible reflecting structures, with the full structure as its last element, as an additional parameter, and then say that the level $V_{\kappa}$ in question which failed to satisfy the reflection principle is ``describable" with respect to this larger set of parameters. Then one might suggest that now there is just as good motivation for taking this as a ``warrant" for going one level further as in the case of the reflection principle discussed in the previous section. This might convince one to accept Roberts' principle as intrinsically justified after all, but on the other hand the well-ordering is now defined in terms of one of the parameters, which might be seen as problematic. However, using the result that every element of $V$ occurs in a set generic extension of $HOD$, we get that every $V_{\kappa}$ does in fact have a canonical well-ordering, and so in this way that objection can be overcome.

\bigskip

To further illuminate the discussion of the problem, it will be useful to analyze further what kinds of comparisons can be made between the reflection principle discussed in the previous section and the reflection principle of Roberts. But, before doing so, we shall introduce the reflection principles discussed by Marshall in \cite{Marshall89} and discuss how they are related to the reflection principle of Roberts (and natural generalizations of it).

\bigskip

So let us now describe Marshall's family of theories $B_n$ indexed by the positive integers $n$, following her exposition in \cite{Marshall89}. The first theory in the sequence, $B_1$, was formulated by Bernays and is presented in \cite{Chuaqui1978}. We shall first describe the theory $B_1$ and then proceed to describe Marshall's generalizations $B_n$ for integers $n>1$.

\bigskip

The theory $B_1$ is a theory in the first-order language of set theory. The objects of the domain of discourse are called classes and a class is said to be a set if and only if it is a member of some class, otherwise it is said to be a proper class. Now we shall describe the axioms of the theory $B_1$.

\bigskip

(A1) Extensionality. $\forall x (x \in A \equiv x \in B) \implies A=B$.

\bigskip

(A2) Class specification. Suppose that $\phi$ is a formula and that $A$ is not free in $\phi$. Then

$$ \exists A \forall x (x \in A \equiv \phi(x) \wedge x \hspace{1 mm} \mathrm{is} \hspace{1 mm} \mathrm{a} \hspace{1 mm} \mathrm{set}). $$

\bigskip

(A3) Subsets. $A \hspace{1 mm} \mathrm{is} \hspace{1 mm} \mathrm{a} \hspace{1 mm} \mathrm{set} \wedge B\subseteq A \implies B \hspace{1 mm} \mathrm{is} \hspace{1 mm} \mathrm{a} \hspace{1 mm} \mathrm{set}.$

\bigskip

(A4) Reflection principle. Suppose that $\phi(x)$ is a formula. Then

$$ \phi(A) \implies \exists u (u \hspace{1 mm} \mathrm{is} \hspace{1 mm} \mathrm{a} \hspace{1 mm} \mathrm{transitive} \hspace{1 mm} \mathrm{set} \wedge \phi^{\mathcal{P}u}(A \cap u)). $$

\bigskip

(A5) Foundation. $\exists x (x \in A) \implies \exists x (x \in A \wedge \forall z (z \notin x \vee z \notin A))$.

\bigskip

(A6) Choice for sets. $\forall x((x \hspace{1 mm} \mathrm{is} \hspace{1 mm} \mathrm{a} \hspace{1 mm} \mathrm{set} \wedge \forall y \forall y' ((y \in x \wedge y' \in x \wedge y \neq y') \implies y\cap y'=\emptyset)) \implies \exists x' (\forall y (y \in x \implies \exists ! z (z \in y \wedge z \in x'))))$.

\bigskip

This completes the description of the axioms of the theory $B_1$. The natural models of $B_1$ are those sets of the form $V_{\alpha+1}$ for some ordinal $\alpha$ such that $(V_{\alpha+1},\in)$ is a model for the theory $B_1$. The ordinals $\alpha$ such that $V_{\alpha+1}$ is a natural model of $B_1$ are precisely the $\Pi^1_\infty$-indescribable cardinals, thus the existence of a $\Pi^1_\infty$-indescribable cardinal is a sufficient assumption on which to prove that $B_1$ is consistent. The theory $B_1$ in turn implies the relativization of all the axioms of $\textsf{ZFC}$ to the class of sets and also the existence of a proper class of $\Pi^1_n$-indescribable cardinals for each positive integer $n$. Let us now see how Marshall generalized this to a sequence of stronger theories. We shall first describe the theory $B_2$ and then proceed to describe the theories $B_n$ for $n>2$.

\bigskip

The theory $B_2$ is a theory in the first-order language of set theory. The objects of the domain of discourse are called 2-classes, an object is said to be a class if and only if it is a member of some 2-class, and an object is said to be a set if and only if it is a member of some class. The class of all sets (whose existence and uniqueness will follow from the axioms soon to be presented) is denoted by $V_0$, and the 2-class of all classes is denoted by $V_1$.

\bigskip

The first axiom for the theory is the axiom of 2-class specification.

\bigskip

(A1) Suppose that $\phi$ is a formula and $A$ is not free in $\phi$. Then

$$ \exists A \forall x (x \in A \equiv \phi(x) \wedge x \hspace{1 mm} \mathrm{is} \hspace{1 mm} \mathrm{a} \hspace{1mm}  \mathrm{class}). $$

Next we have extensionality for 2-classes.

\bigskip

(A2) $\forall x (x \in A \equiv x \in B) \implies A=B$.

\bigskip

We write $\mathrm{STC}(v)$ to abbreviate a formula saying that $v$ is transitive, and every sub-2-class $x$ of any element $y$ of $v$ is also an element of $v$.

\bigskip

If $A$ and $u$ are 2-classes, $A^u$ denotes $A \cap u$ if $A$ is a class, and denotes $\{x^u\mid x \in A \cap u\}$ if $A$ is not a class. $\phi^u$ denotes the relativization of the formula $\phi$ to $u$.

\bigskip

If $x$ and $y$ are two classes, then the ordered pair $[x,y]$ of the classes $x$ and $y$ is defined to be $x \times \{0\} \cup y \times \{1\}$. The notation $A(u)$ is an abbreviation for $u \cap V_0 \in V_0 \wedge \mathrm{STC}(u \cap V_0) \wedge \forall x \forall y (x, y \in u \implies [x,y] \in u)$. Then our reflection principle is given by

\bigskip

(A3) If $\phi(x)$ is a formula, then

$$ \phi(A) \implies \exists u (A(u) \wedge \phi^{\mathcal{PP}(u \cap V_0)}(A^u)). $$

Finally we take as axioms the axiom of foundation and the axiom of choice for sets as in the theory $B_1$, noting that the interpretation of the abbreviation ``$x$ is a set" is now different. This completes the description of the theory $B_2$. Given any theorem $\phi$ of $B_1$, the relativization of $\phi$ to $V_1$ is a theorem of $B_2$, as shown in \cite{Marshall89}. Marshall also shows in \cite{Marshall89} that $B_2$ proves that the class of 1-extendible cardinals is a stationary subclass of the class of ordinals, and that if $\kappa$ is 2-extendible or $\beth_{\kappa+1}$-supercompact, then $V_{\kappa+2}$ is a natural model of $B_2$. Before proceeding to give the description of the theories $B_n$ with $n>2$, let us pause to discuss the relationship between $B_2$ and a natural generalization of Roberts' theory $\textsf{ZFC2}_{S}+\textsf{R}_{S}$, and also to characterize the natural models of $B_2$ and related theories, thereby answering a question asked by Marshall in \cite{Marshall89}.

\bigskip

There is an obvious extension of Roberts' theory $\textsf{ZFC2}_{S}+\textsf{R}_{S}$ to a theory in the third-order language of set theory $\textsf{ZFC3}+\textsf{R}_{3}$. And there is an obvious interpretation of the third-order language of set theory in the first-order language of set theory where the domain of discourse of the variables in the latter language is understood to be the 2-classes. It is easy to see that the image of the theory $\textsf{ZFC3}+\textsf{R}_{3}$ under this interpretation is equivalent to the theory $B_2$. Similar relationships could be found between other natural extensions of Roberts' theory and the theories $B_n$ for integers $n>2$ to be defined later. Let us characterize the natural models of $B_2$. Roberts shows in \cite{Roberts2017} that if $V_{\kappa+2}$ is a natural model of $B_2$ then $\kappa$ is a limit of a 1-extendible chain, that is, a subset $S \subseteq \kappa$ cofinal in $\kappa$ such that if $\alpha \in S, \alpha<\beta$, and $\beta \in S \cup \{\kappa\}$, then $\alpha$ is 1-extendible to $\beta$, that is, there is an elementary embedding $j: V_{\alpha+1} \prec V_{\beta+1}$ with critical point $\alpha$ and $j(\alpha)=\beta$. Given any positive integer $n$ there is a clear corresponding notion of a $\Pi^2_n$-extendible chain, and Roberts' result easily generalizes to show that $V_{\kappa+2}$ is a natural model of $B_2$ if and only if $\kappa$ is a limit of a $\Pi^2_n$-extendible chain for all positive integers $n$. There is a similar characterization of the natural models of $\textsf{ZFC2}_{S}+R_{S}$ in which $V_{\alpha+2}$ is replaced with $L_{1}(V_{\alpha+1})$ and one uses the hierarchy of complexity classes of formulas relativized to $L_{1}(V_{\alpha+1})$. Let us now return to the task of defining the theories $B_n$ for integers $n>2$.

\bigskip

The theory $B_n$ for an integer $n>2$ is a theory in the first-order language of set theory where the domain of discourse is taken to be the $n$-classes, where classes are collections of sets, 2-classes are collections of classes, 3-classes are collections of 2-classes, and so on. The axioms of extensionality and class specification are as before, and there is an obvious definition of $V_0$, the collection of sets, $V_1$, the collection of classes, and $V_k$, the collection of $k$-classes for an integer $k$ such that $2 \leq k \leq n$. For arbitrary $n$-classes $A$ and $u$, we define $A^u=A \cap u$ if $A \in V_1$ and then by $\in$-induction we define $A^u=\{x^u\mid x \in A \cap u\}$ for $A \notin V_1$. Then our reflection principle is: For each formula $\phi(x)$, we have

$$ \phi(A) \implies \exists u (A(u) \wedge \phi^{\mathcal{P}^n(V_0 \cap u)}(A^u)). $$

Then we have the axiom of foundation and the axiom of choice for sets as before. This completes the definition of the theory $B_n$, which is equivalent to natural extensions of Roberts' theory $\textsf{ZFC2}_{S}+\textsf{R}_{S}$ in the same way as before.

\bigskip

Suppose that we consider an extension of Roberts' theory to a theory in the $\omega$th-order language of set theory which we might label $\textsf{ZFC}\omega+\textsf{R}_{\omega}$, and suppose that we drop the requirement that the reflecting structure must be ``set-sized", but keep the requirement that it must be ``standard". The natural models of this theory are of the form $V_{\kappa+\omega}$ where $\kappa$ is a limit of an $n$-extendible chain for each positive integer $n$. We could also easily formulate an essentially equivalent theory in the first-order language of set theory along the lines outlined by Marshall, where now the role of what corresponds to what we previously called the ``reflecting structure" is no longer to serve as the reflecting structure but simply to guide the reflection of the parameter.

\bigskip

On the other hand, I might modify the reflection principle in the $\omega$th-order language of set theory described in Section 2, so that the reflecting structure is allowed to be any structure of the form $\langle V_{\alpha}, X_1, X_2, \ldots \rangle$ where $\alpha<\kappa$ and $X_i \subseteq V_{\kappa+i}$ for $i\in\omega$ and $i>0$, and each $X_i$ is standard for $V_{\alpha}$ in the sense that the image under the transitive collapsing map contains all of $V_{\alpha+i}$. That is, the permitted reflecting structures are exactly the same as those allowed in the previously described modification of Roberts' theory. The resulting theory is in fact equivalent to the theory originally described in Section 2, and the natural models of my theory are still the $\omega$-reflective cardinals. The only difference between my theory and the modification of Roberts' theory is that my theory does not allow reflection of formulas which have unbounded existential quantification of first-order variables within the scope of a quantifier for a higher-order variable. However, dropping that requirement on the formulas is only consistent if we allow a wide range of possible reflecting structures.

\bigskip

This indicates the kind of comparison that can be made between the reflection principles that give rise to $\omega$-reflective cardinals, and the much stronger reflection principles considered by Roberts and Marshall which give rise to $n$-extendible cardinals for every positive integer $n$. In \cite{McCallum2017} we showed how the property of being a remarkable cardinal could be seen as a reflection principle. As discussed in more detail below, if we construct a reflection principle which is to the reflection principle for the remarkable cardinals as this just mentioned stronger reflection principle is to the $\omega$-reflective cardinals, then we obtain a large-cardinal property equivalent to supercompactness.

\bigskip

In this way we can see two different possible notions of intrinsic justification, one extending all the way up to a remarkable cardinal, and one extending all the way up to a supercompact cardinal.

\bigskip

As mentioned in the discussion to which I referred just now, of the connection between the notion of a remarkable cardinal and the notion of a supercompact cardinal, it is natural to generalise, along a line of thought similar to that which yielded the generalization from $\omega$-reflective cardinals to the reflection principle for remarkable cardinals, from Roberts' reflection principle to a stronger reflection principle which yields the existence of supercompact cardinals. Given a level $V_{\kappa}$, we can consider structures of the form $(V_{\kappa}, V_{\lambda})$ with $\lambda>\kappa$ and consider some formula $\phi$ in a two-sorted language holding in such a structure relative to a certain finite collection of parameters. It is natural to posit that there should exist a ``set-sized" reflecting structure, containing all the parameters, whose first component is $V_{\alpha}$ for some $\alpha<\kappa$ and whose second component is ``set-sized" in the sense of having cardinality less than $\beth_{\kappa}$ and furthermore such that the transitive collapse of the second component is of the form $V_{\beta}$ for some $\beta>\alpha$. (Here the collapsing map may not be injective.) A level $V_{\kappa}$ satisfies this form of reflection if and only if $\kappa$ is supercompact. This form of reflection is quite closely related to, but not quite as strong as, the form of reflection considered in Marshall's theory $B_0(V_0)$ of the paper \cite{Marshall89}. Indeed, we easily see from the foregoing considerations that Marshall's theory $B_0(V_0)$ proves that the height $\kappa$ of the universe $V_0$ is supercompact, and therefore by further application of the reflection principle that $V_0$ is a model for the assertion that there is a proper class of supercompact cardinals. This resolves affirmatively the question asked by Marshall about whether $B_0(V_0)$ is sufficiently strong to prove the existence of a supercompact cardinal. On the other hand it would seem that her theory $B_1(V_0)$ does not prove the existence of a huge cardinal. For if $\kappa$ is huge and $j:V\prec M$ is an elementary embedding witnessing the hugeness of $\kappa$, then it seems that $(V_{\kappa}, V_{j(\kappa)}, V_{j(\kappa)+1})$ is a model for $B_1(V_0)$ using the same reasoning which Marshall uses to prove that it is a model of $B_0(V_0)$, in the proof Theorem 3 of Chapter IV of \cite{Marshall89}.

\bigskip

Later in the paper \cite{Marshall89} Marshall introduces the theories $B_n(V^0_0, V^1_0, \ldots V^{m-1}_{0})$, for positive integers $m$ and $n$, as motivation for the existence of $m$-huge cardinals. In the final section of her paper she goes even further up to the point of inconsistency with the axiom of choice. I hope in later work to elaborate on possible principled reasons, inspired by the lines of thought put forward by Marshall, and assuming that one does indeed regard the principles of Welch and Roberts as intrinsically justified, for accepting all large cardinals not known to be inconsistent with the axiom of choice, while still finding principled reasons for not going to the point of inconsistency with choice.

\section{Concluding remarks}

Our exploration of the range of reflection principles which have been proposed in the literature outline two possible conceptions of intrinsic justification, one which matches Peter Koellner's conjecture that intrinsic justifications do not take us to the point of an $\omega$-Erd\H{o}s cardinal, and one which extends up to the reflection principles of Welch and Roberts, and indeed, as suggested here, plausibly all the way up to a supercompact cardinal. The remarks made here should serve as a basis for further philosophical exploration of the question of which large cardinals are intrinsically justified.

\pagebreak[4]

\end{document}